\documentclass[12pt]{amsart}
\usepackage{epic}
\usepackage[all]{xy}
\usepackage{amsfonts}
\usepackage{amssymb}

\usepackage{color}

\newtheorem{teo}{Theorem}[section]
\newtheorem{lem}[teo]{Lemma}

\newtheorem{cor}[teo]{Corollary}

\newtheorem{assum}[teo]{Assumption}
\theoremstyle{definition}
\newtheorem{dfn}[teo]{Definition}
\newtheorem{rk}[teo]{Remark}

\def\<{\langle}
\def\>{\rangle}
\def\ss{\subset}

\def\a{\alpha}
\def\b{\beta}
\def\g{\gamma}

\def\e{\varepsilon}

\def\t{\tau}

\def\Z{{\mathbb Z}}

\def\Aut{\operatorname{Aut}}

\def\1{\mathbf 1}

\newcommand{\til}[1]{\widetilde{#1}}
\newcommand{\wh}[1]{\widehat{#1}}

\def\cT{{\mathcal T}}

\def\N{{\mathbb N}}

\def\St{\operatorname{St}}
\def\Iso{\operatorname{Iso}}
\def\WBI{\operatorname{WBI}}

\begin{document}

\title[Reidemeister numbers of weakly branch groups]
{Reidemeister numbers of saturated weakly branch groups}

\author{Alexander Fel'shtyn}
\address{Instytut Matematyki, Uniwersytet Szczecinski,
ul. Wielkopolska 15, 70-451 Szczecin, Poland and Department of Mathematics,
Boise State University,
1910 University Drive, Boise, Idaho, 83725-155, USA}
\email{felshtyn@diamond.boisestate.edu, felshtyn@mpim-bonn.mpg.de}

\author{Yuriy Leonov}
\address{IT Department, Odessa Academy of Telecommunications,
Kuznechnaja street, 1, 65000 Odessa, Ukraine}
\email{leonov\_yu@yahoo.com}

\author{Evgenij Troitsky}
\thanks{The present research is a part of joint research programm in Max-Planck-Institut
f\"ur Mathematik (MPI) in Bonn with partial relation to the activity
``Geometry and Group Theory''. We would like to thank the MPI for its kind support and
hospitality while the most part of this work has been completed.\\
\hspace*{0.5em} The third author is partially supported by
RFFI Grant  05-01-00923
and Grant ``Universities of Russia''}
\address{Dept. of Mech. and Math., Moscow State University,
119992 GSP-2  Moscow, Russia}
\email{troitsky@mech.math.msu.su}
\urladdr{
http://mech.math.msu.su/\~{}troitsky}


\begin{abstract}
We prove for a wide class of saturated weakly branch group
(including the (first) Grigorchuk group and the Gupta-Sidki group)
that the Reidemeister
number of any automorphism is infinite.
\end{abstract}

\maketitle

\tableofcontents

\section{Introduction}
Let $\phi:G\to G$ be an automorphism of a group $G$. A class of equivalence
$x\sim gx\phi(g^{-1})$ is called the \emph{Reidemeister class }
or $\phi$-\emph{conjugacy class} or \emph{twisted conjugacy class of} $\phi$. The
number $R(\phi)$ of Reidemeister classes is called the \emph{Reidemeister number.}
The interest in twisted conjugacy relations has its origins, in particular,
in the Nielsen-Reidemeister fixed point theory (see, e.g. \cite{Jiang,FelshB}),
in Selberg theory (see, eg. \cite{Shokra,Arthur}), and  Algebraic Geometry
(see, e.g. \cite{Groth}). The main current problem of the field is to obtain
a twisted analogue of the celebrated Burnside-Frobenius theorem
\cite{FelHill,FelshB,FelTro,FelTroVer,ncrmkwb,polyc,FelTroObzo}. For this purpose
it is important to describe the class of groups $G$, such that $R(\phi)=\infty$ for
any automorphism $\phi:G\to G$. First attempts to localize this class of groups go up to
\cite{FelHill}.
After that
it was proved that the following groups belong to this class:
non-elementary Gromov hyperbolic groups \cite{FelPOMI,ll},
Baumslag-Solitar groups $BS(m,n) = \langle a,b | ba^mb^{-1} = a^n \rangle$
except for $BS(1,1)$ \cite{FelGon},
generalized Baumslag-Solitar groups, that is, finitely generated groups
which act on a tree with all edge and vertex stabilizers infinite cyclic
\cite{LevittBaums},
lamplighter groups $\Z_n \wr \Z$ iff $2|n$ or $3|n$ \cite{gowon1},
the solvable generalization $\Gamma$ of $BS(1,n)$ given by the short exact sequence
$$1 \rightarrow \Z[\frac{1}{n}] \rightarrow \Gamma \rightarrow \Z^k \rightarrow 1$$
 as well as any group quasi-isometric to $\Gamma$ \cite{TabWong},
groups which are quasi-isometric to $BS(1,n)$ \cite{TabWong2} (while this property is
 not a quasi-isometry invariant), the chameleon
 R.~Thompson group \cite{ThomsonRinf}.

In paper \cite{TabWong} a terminology for this property was suggested.
Namely, a group $G$ has \emph{property } $R_\infty$ if any its automorphism $\phi$
has $R(\phi)=\infty$.

For the immediate consequences of $R_\infty$ property for the topological fixed point
theory e.g. see \cite{TabWong2}.

In the present paper we prove that a wide class of weakly branch groups including
the Grigorchuk group and the Gupta-Sidki group, has $R_\infty$ property.

The results of the present paper demonstrate that the further study
of Reidemeister theory for this class of groups has to go
along the lines specific for the infinite case. On the other
hand these results make smaller the class of groups, for which the
twisted Burnside-Frobenius conjecture \cite{FelHill,FelTro,FelTroVer,ncrmkwb,polyc,FelTroObzo}
has to be verified.

\medskip\noindent
{\bf Acknowledgement.}
The authors are grateful to
R.~Grigorchuk,
Yu.~Prokhorov,
M.~Sapir,
S.~Sidki and
Z.~\v Sunic
for helpful discussions.

\section{Preliminaries on weakly branch groups}
Let $\cT$ be a (spherically symmetric) rooted tree.
A group $G$ acting faithfully on a rooted tree is said to be
a \emph{weakly branch group} if for every vertex $v$ of $\cT$
there exists an element of $G$ which acts nontrivially on the
subtree $\cT_v$ with the root vertex $v$ and trivially off it.

The group $G$ is said to be \emph{saturated } if for every
positive integer $n$ there exists a characteristic subgroup
$H_n \ss G$ acting trivially on the $n$-th level of $\cT$ and
level transitive on any subtree $\cT_v$ with $v$ in the $n$-th
level.

\begin{teo}[\cite{LavrNekr}]\label{teo:lavnik}
Let $G$ be a saturated weakly branch group. Then the
automorphism group $\Aut G$ coincides with the normalizer of $G$
in the full automorphism group $\Iso\cT$ of the rooted tree; i.e., every
automorphism of the group $G$ is induced by conjugation from
the normalizer and the centralizer of $G$ in $\Iso\cT$ is trivial.
\end{teo}

All groups in what follows will be supposed to be
\emph{saturated weakly branch.}


\begin{teo}[{\cite[Theorem 7.1]{LavrNekr}}]\label{teo:lavnek.7.1}
Let $G$ be a level-transitive isometry group of a rooted tree $\cT$
with a stabilizer sequence $(G_0\supseteq G_1 \supseteq G_2\supseteq\dots )$.
An
automorphism $\phi\in \Aut G$ is induced by an element of
the full automorphism group $\Iso\cT$ of the rooted tree
if and only if for every $n\ge 0$ there exists $a_n\in G$ such that
$a_n G_i a_n^{-1}=\phi(G_i)$ for every $i\le n$.
\end{teo}

\section{Reidemeister classes and inner automorphisms}

Let us denote by $\t_g:G\to G$ the automorphism $\t_g(\til g)=g\til g\,g^{-1}$
for $g\in G$. Its restriction on a normal subgroup we will denote by $\t_g$
as well. We will need the following statements.

\begin{lem}\label{lem:redklassed}
$\{g\}_\phi k=\{g\,k\}_{\t_{k^{-1}}\circ\phi}$.
\end{lem}

\begin{proof}
Let $g'=f\,g\,\phi(f^{-1})$ be $\phi$-conjugate to $g$. Then
$$
g'\,k=f\,g\,\phi(f^{-1})\,k=f\,g\,k\,k^{-1}\,\phi(f^{-1})\,k
=f\,(g\,k)\,(\t_{k^{-1}}\circ\phi)(f^{-1}).
$$
Conversely, if $g'$ is $\t_{k^{-1}}\circ\phi$-conjugate to $g$, then
$$
g'\,k^{-1}=
f\,g\,(\t_{k^{-1}}\circ\phi)(f^{-1})k^{-1}=
f\,g\,k^{-1}\,\phi(f^{-1}).
$$
Hence a shift maps $\phi$-conjugacy classes onto classes related to
another automorphism.
\end{proof}

\begin{cor}\label{lem:innerreidem}
$R(\phi)=R(\t_g \circ \phi)$.
\end{cor}

\section{Binary trees}
Consider the subset $K_n \ss  G$ formed of elements
of the stabilizer $\St_n G$ such that on the $n+1$-th level they are switching
each pair of (neighboring) vertexes.

\begin{lem}\label{lem:knisnonempty}
Under our assumptions, $K_n$ is not empty.
\end{lem}

\begin{proof}
Let us encode the action of an element of $G$, which consists only of switches
or trivial actions on neighboring vertices as a sequence of $-1$ and $1$.
For $n=1$ the statement is evident. Let us argue by induction and suppose
that the statement is true for the levels up to $n-1$.
By the supposition on $G$ to be saturated, among the elements of $\St_n(G)$ there
are some elements, such that for any pair of neighboring vertices one of
these elements has $-1$ on the corresponding places. We have two possibilities:
in any pair of neighbors both entries are equal to each other, or there is a pair
with $+1$ and $-1$. In the second case let us conjugate our element by the desired
element at the level $n-1$. The result of its action is the following: it transposes
each pair. Indeed, since the automorphism group of $\Z_2$ is $\Z_2$, the conjugation
can only permutate. For example, the conjugation sends
$$
((1,-1),(1,1),(-1,-1),(1,-1),\dots)
$$
to
$$
((-1,1),(1,1),(-1,-1),(-1,1),\dots).
$$
Their product will have $(-1,-1)$ on the place under consideration and pairs of the
same elements (i.e., $(-1,-1)$ or $(1,1)$) on the remaining ones.
So, we have reduced the second possibility to the first one.
We go further
taking the conjugation by the element, which was obtained at the level $n-2$.
After an analogous multiplication we obtain
a (nontrivial) element with quadruples of neighbors formed by the same
elements. And so
on. The end step of the induction (may be the first one) is the desired element at
the level $n$.
\end{proof}

\begin{lem}\label{lem:mainlemfortt}
Let an automorphism $t$ $($after the elimination of action on upper levels$)$
have at level $m+1$ the number of switches, which is distinct from $2^{m-1}$.
Then the Reidemeister class of an element from $K_m$ does not intersect
$\St_{m+1}G$.
\end{lem}

\begin{proof}
Suppose, $g\in K_m$, $h\in G$. Let us consider
\begin{equation}\label{eq:twistconjelem}
hg\phi(h^{-1})=hg th^{-1} t^{-1}=(hgh^{-1})(hth^{-1} t^{-1}).
\end{equation}
By the same argument with $\Z_2$ as in proof of Lemma \ref{lem:knisnonempty},
$hgh^{-1}\in K_m$. If $hth^{-1} t^{-1}$ is non-identical on some of levels $1,\dots,m$,
then the product (\ref{eq:twistconjelem}) is non-identical as well, since
$hgh^{-1}\in K_m$, and we are done.
Otherwise, let us remark, that after elimination of the action on the previous levels
the number of switches on the level $m+1$ in $hgth^{-1}$ is $2^m$-(the number of switches
of $t$), which is distinct from the number of switches of $t$ (or $t^{-1}$)
provided that it is not equal to $2^{m-1}$. Hence the total number of switches in
 (\ref{eq:twistconjelem}) is non-zero. So it is non-trivial on the level $m+1$.
\end{proof}

\begin{teo}\label{teo:reidinftybinary}
Let $G$ be a saturated weakly branch group acting on a binary tree $\cT$.
Suppose, $\phi:G\to G$ is an automorphism, such that for any $k\in \N$
there exists an inner automorphism of $G$ such that its composition $\phi'$ with
$\phi$ satisfies the condition of Lemma {\rm\ref{lem:mainlemfortt}} at some
collection of levels of number $k$. Then $R(\phi)=\infty$.
\end{teo}

\begin{proof}
Let us take an arbitrary $k\in\N$ and show that the number of Reidemeister
classes is not less than $k$. For this purpose, take an appropriate
inner automorphism in accordance with the supposition. Then the Reidemeister
 numbers of $\phi$ and $\phi'$ are the same (cf. Lemma
\ref{lem:innerreidem}). So it is sufficient to prove that $R(\phi')\ge k$.
For the notation brevity, suppose that the levels with mentioned parity properties
are $1,\dots,k$.
By Lemma \ref{lem:mainlemfortt}, for elements $g_i\in K_i\ss \St_i G$
one has
$\{g_i\}_{\phi'}\cap \St_{i+1} G=\varnothing$, $i=1,\dots,k$. Hence
the classes $\{g_i\}_{\phi'}$, $i=1,\dots,k$,
are distinct.
\end{proof}

\begin{rk}
It is clear, that the following condition can serve as
an alternative for the supposition of the theorem: there exists a
$g_i\in K_i$ with an odd number of switches.
\end{rk}

\section{Grigorchuk group}\label{sec:Griggroup}

Now we want to prove that the Grigorchuk group $G$ (\cite{GrFA},
see also \cite{harpe,GrigorHorizon,GrigSidki}) satisfies the conditions
of Theorem \ref{teo:reidinftybinary}.

Consider the following presentation of $G$ (cf. \cite{GrigSidki}). It
has the generators $a,b,c$, and element $d$, where $d=bc$, $a$ is defined
at Fig.~\ref{fig:defgriga},
where $-1$ is a switch,
$b$ and $c$ are defined inductively by
$$
b=(a,c),\qquad c=(a,d),\qquad d=(1,b),
$$
where brackets mean the action on the corresponding sub-trees.
\begin{figure}[hb]
\unitlength=1mm
\begin{picture}(22,18)(0,0)
\thicklines
\put(2,5){\line(1,1){9}}
\put(22,5){\line(-1,1){9}}
\put(0,0){$1$}
\put(20,0){$-1$}
\put(11,15){$1$}
\end{picture}
\caption{}\label{fig:defgriga}
\end{figure}

Then $b$ and $c$ are as at Fig.~\ref{fig:defgrigbc}  (we partially omit $1$'s).
\begin{figure}[ht]
\unitlength=1mm
\begin{picture}(100,35)(0,0)
\thicklines
\put(5,35){\line(-1,-1){5}}
\put(10,30){\line(-1,-1){5}}
\put(15,25){\line(-1,-1){5}}
\put(20,20){\line(-1,-1){5}}
\put(25,15){\line(-1,-1){5}}
\put(30,10){\line(-1,-1){5}}
\put(35,5){\line(-1,1){30}}
\dottedline[.]{1.0}(36,4)(39,1)
\dottedline[.]{1.0}(34,4)(31,1)
\put(23,3){$\scriptstyle 1$}
\put(16,7){$\scriptstyle -1$}
\put(11,12){$\scriptstyle -1$}
\put(8,18){$\scriptstyle 1$}
\put(1,22){$\scriptstyle -1$}
\put(-4,27){$\scriptstyle -1$}
\put(55,35){\line(-1,-1){5}}
\put(60,30){\line(-1,-1){5}}
\put(65,25){\line(-1,-1){5}}
\put(70,20){\line(-1,-1){5}}
\put(75,15){\line(-1,-1){5}}
\put(80,10){\line(-1,-1){5}}
\put(85,5){\line(-1,1){30}}
\dottedline[.]{1.0}(86,4)(89,1)
\dottedline[.]{1.0}(84,4)(81,1)
\put(71,3){$\scriptstyle -1$}
\put(68,7){$\scriptstyle 1$}
\put(61,12){$\scriptstyle -1$}
\put(56,18){$\scriptstyle -1$}
\put(53,22){$\scriptstyle 1$}
\put(46,27){$\scriptstyle -1$}
\put(0,15){$\scriptstyle b$}
\put(50,15){$\scriptstyle c$}
\end{picture}
\caption{}\label{fig:defgrigbc}
\end{figure}

In particular
\begin{equation}\label{eq:relatGrig}
    a^2=b^2=c^2=d^2=1,\qquad a^{-1}c a =(d,a),\qquad a^{-1}d a =(b,1).
\end{equation}

By \cite{GrigSidki}, any automorphism of Grigorchuk group, up to
taking a product by an inner one, is a finite product of
commuting involutions of the form (at some level):
$$
(1, (ad)^2,1, (ad)^2,1, (ad)^2,1, (ad)^2,\dots).
$$
Hence, by (\ref{eq:relatGrig}) on the next level we will have:
$$
(1,1,b,b,1,1,b,b,1,1,b,b,\dots).
$$
Keeping in mind the form of $b$ (see Fig.~\ref{fig:defgrigbc}), we conclude
that the number of switches at low levels is bounded (by the number of $b$'s
in the above formula). Hence, some uniform estimation holds for their
finite product (i.e., our automorphism). Thus, starting from some level, the
number of switches is less than the half of places, and we can apply
Theorem \ref{teo:reidinftybinary}.

\section{Some generalizations: strongly saturated groups}
We return to the case of a general spherical tree $\cT$.

\begin{dfn}
Let us denote by $l(m)$ the number of vertexes at the level $m$.
\end{dfn}

We will make the following supposition about $t$.
\begin{assum}\label{assum:goodelem}
There exist a constant $s\in (0,1)$ such that
for any $j$ the isometry $t$ has the
number of fixed vertexes at the level $j$ not less than $s\cdot l(j)$.
\end{assum}

We will need the following definition.

\begin{dfn}\label{dfn:strongsatur}
A saturated group $G$ is called \emph{strongly saturated}
if for any $i\in \N$ there exists an element $g_i\in \St_i(G)$
such that it has no fixed points on the level $i+1$.
\end{dfn}

As it is proved above, any saturated group on a binary tree
is strongly saturated. We will see (Remark \ref{rk:GupSidstrsat})
that the Gupta-Sidki group is strongly saturated as well.

\begin{lem}\label{lem:vlastelinkolez}
Suppose, an automorphism $\phi$
of a strongly saturated group $G$
is defined by $t\in\Iso\cT$ and
Assumption {\rm\ref{assum:goodelem}} holds.
Then there exists an element $\wh g\in \St_m G$ such that its Reidemeister
class $\{\wh g\}_\phi$ does not intersect $\St_{m+r}(G)$, where
$$
\left(\frac 12\right)^r< s.
$$
\end{lem}

\begin{proof}
An element of $\{\wh g\}_\phi$ has the form $h\wh g t h^{-1} t^{-1}=:\wh g_h$. Note
that for any $\wh g\in \St_m G$ this element is not in $\St_{m}G$ and hence
not in $\St_{m+l(m)}(G)$, if $h th^{-1}t^{-1}\not\in\St_{m}G$. Indeed,
$$
\wh g_h=(h\wh g h^{-1}) (ht h^{-1} t^{-1}),
$$
where the first factor belongs $\St_m G$. Thus we need to check for the
desired $\wh g$ only that
\begin{equation}\label{eq:propghat}
    \wh g_h=(h(\wh g t) h^{-1}) t^{-1}\not\in\St_{m+r}(G),\qquad \mbox{ for $h$ such that }
h th^{-1}t^{-1}\in\St_{m}G,
\end{equation}
i.e., $h$ and $t$ commute on the level $m$. In this case one has
$$
t(h v_m)=h t v_m= h v_m,
$$
i.e. $h$ maps $v_m$ to some other fixed point of $t$ (cf. Fig.~\ref{fig:towers1},
where thick points mean some fixed points of $t$).
\begin{figure}[ht]
\unitlength=1bp
\begin{picture}(245,135)(0,0)
\put(0,0){
}
\put(70,0){
}
\put(150,0){
}
\put(0,120){\line(1,0){220}}
\put(-8,0){\line(1,0){245}}
\put(222,118){$\scriptstyle m$}
\put(240,-4){$\scriptstyle m+l(m)$}
\put(30,120){\circle*{5}}
\put(100,120){\circle*{5}}
\put(180,120){\circle*{5}}
\put(182,112){\circle*{5}}
\put(184,104){\circle*{5}}
\put(186,96){\circle*{5}}
\put(208,8){\circle*{5}}
\put(210,0){\circle*{5}}
\put(179,125){$\scriptstyle v_m$}
\put(29,125){$\scriptstyle v''$}
\put(99,125){$\scriptstyle v'$}
\put(186,112){$\scriptstyle v_{m+1}$}
\put(188,104){$\scriptstyle v_{m+2}$}
\put(190,96){$\scriptstyle v_{m+3}$}
\put(209,11){${\scriptstyle v_{m+l(m)-1}}$}
\put(211,4){${\scriptstyle v_{m+l(m)}}$}
\put(140,131){$\scriptstyle h$}
\put(140,105){$\scriptstyle h^{-1}$}
\bezier{300}(100,120)(140,110)(180,120)
\bezier{300}(100,120)(140,130)(180,120)
\put(140,123){$\scriptstyle \blacktriangleleft$}
\put(140,112){$\scriptstyle \blacktriangleright$}
\end{picture}
\caption{}\label{fig:towers1}
\end{figure}

We will construct the desired $\wh g$ as a product
$$
\wh g= g_{m+l(m)}^{\e(r)}  \cdots g_{m+1}^{\e(1)}  g_{m}^{\e(0)},\quad
\mbox{ where $\e(i)=0$ or $1$, } x^0:=e,
$$
where $g_i$ are from Definition \ref{dfn:strongsatur}.
We will chose these $\e(i)$ in such a way that $\wh g t$ will have the number of
fixed points in the set $B_r$ of the vertexes of
the level $m+r$ is less than $s\cdot l(m+r)=s\cdot \# B_r$.
This $\wh g$ will be
the desired one. Indeed, the map $t$ on the $m+r$ level
has at least $s\cdot l(m+r)$
fixed points, and to obtain the identity in the composition
(\ref{eq:propghat}) one should suppose that $\wh g t$ has the same number of fixed points
at the level $m+r$.

Now we pass to the determination of $\e(i)$ (cf. Fig.~\ref{fig:bashnia}).
If the number of fixed points
of $t$ at the set $B_1$ of all vertexes of $\cT$ at the level $m+1$ is less
than 1/2 of $\# B_1$, then we take $\e(0)=0$. Otherwise, $\e(0)=1$. Then
$g_{m}^{\e(0)} t$ has the number of fixed points in $B_1$ less or equal than
$\frac 12 \# B_1$. Define $B_2$ in the same way as $B_1$, and let $C_2\ss B_2$
be the set of that points, which do not belong to subtrees with roots
at level $m+1$, which are not fixed by $g_{m}^{\e(0)} t$. Thus,
$\# C_2 \le \frac 12 \# B_2$. The points of $B_2\setminus C_2$ are
not-fixed by $g_{m}^{\e(0)} t$. Now we consider the fixed points of
 $g_{m}^{\e(0)} t$ in $C_2$. If the number of them is less than $\frac 12\# C_2$,
 then we take $\e(1)=0$. Otherwise, $\e(1)=1$.
 Take the composition $g_{m+1}^{\e(1)}  g_{m}^{\e(0)} t$. Since $g_{m+1}\in\St_{m+1}G$,
 the vertexes which come from the vertexes at the level $m+1$, which were not
 fixed by $g_{m}^{\e(0)} t$, are still not fixed by $g_{m+1}^{\e(1)}  g_{m}^{\e(0)} t$.
Hence, if $C_3\ss B_3$ is defined as the complement to the set formed
by the vertexes of that subtrees, whose roots were
moved by $g_{m}^{\e(0)} t$ at the level $m+1$ and by $g_{m+1}^{\e(1)}  g_{m}^{\e(0)} t$
at the level $m+2$, then
\begin{enumerate}
    \item $\# C_3\le \left(1-\frac 12 -\frac 14\right) \# B_3$;
    \item $g_{m+1}^{\e(1)}  g_{m}^{\e(0)} t$ has no fixed points on $B_3\setminus C_3$.
\end{enumerate}
\begin{figure}[ht]
\unitlength=1bp
\begin{picture}(180,180)(0,0)
\put(0,0){
}
\put(0,0){
}
\put(90,0){
}
\put(90,0){
}
\put(90,0){
}
\put(90,0){
}
\put(90,0){
}
\put(90,150){\line(1,0){10}}
\put(90,120){\line(1,0){10}}
\put(90,90){\line(1,0){11}}
\put(100,60){\line(1,0){5}}
\put(93,154){$\scriptstyle C_2$}
\put(91,124){$\scriptstyle C_3$}
\put(137,104){$\scriptstyle C_4$}
\put(155,74){$\scriptstyle C_5$}
\put(170,44){$\scriptstyle C_6$}
\put(135,106){\line(-1,0){40}}
\put(153,76){\line(-1,0){54}}
\put(168,46){\line(-1,0){63.5}}
\put(99,76){\circle*{3}}
\put(104,46){\circle*{3}}
\put(95,106){\circle*{3}}
\end{picture}
\caption{}\label{fig:bashnia}
\end{figure}

Now we count the number of fixed points of $g_{m+1}^{\e(1)}  g_{m}^{\e(0)} t$ on $C_3$, and
so one. Since at each step at least a half of points from $C_i$ comes to ``the world of
non-fixed points'', while the points which are in the subtrees, whose roots were
``joined to the world of non-fixed points'' at the previous steps, can not ``leave this world'',
because $g_i \in \St_{i-1}(G)$, we will obtain the desired result in (no more than) $r$
steps. Indeed, the number of fixed points of $\wh g t$ on the level $m+r$ is less than
$$
\# C_r \le \left(\frac 12\right)^r \# B_r = \left(\frac 12\right)^r l(m+r),
$$
while, by the supposition,  the number of fixed points of $t$ on the level
$m+r$ is more than
$$
s\cdot l(m+r)> \left(\frac 12\right)^r l(m+r).
$$
\end{proof}

\begin{teo}\label{teo:maingenerbigt}
Suppose, $G$ is a weakly branch group on a spherical tree $\cT$
and $\phi$ its automorphism induced by $t\in \Iso\cT$
restricted to satisfy Assumption {\rm \ref{assum:goodelem}}.
Then $R(\phi)=\infty$.
\end{teo}

\begin{proof}
Let us prove that $R(\phi)\ge n$ for any $n\in \N$.

Now we apply Lemma \ref{lem:vlastelinkolez} for the purpose to produce
elements $\wh g_0,\dots,\wh g_{n-1}$, such that
\begin{enumerate}
    \item $\wh g_i\in \St_{ir}(G)$, $i=0,\dots,n-1$;
    \item $\{\wh g_i\}_{\phi'} \cap \St_{(i+1)r}(G)=\varnothing$, $i=0,\dots,n-1$.
\end{enumerate}
Hence, $R(\phi')\ge n$.
\end{proof}

\begin{rk}
In fact, we need much weaker assumptions, than \ref{assum:goodelem}. For example,
one can suppose the existence of
a large number of fixed vertexes, say at each $k$-th level., etc.
This is the case of the Grigorchuk group, as it is evident from Section
\ref{sec:Griggroup}.

Also, it is possible to suppose that Assumption \ref{assum:goodelem}
holds for $gt$ for some $g$, etc.
\end{rk}

\section{Further results: locally normal groups}
We will introduce the following definition, which is related to the property
of $G$ to be residually finite (cf. \cite[p.~171]{LavrNekr}).

\begin{dfn}\label{dfn:normality}
Let $G$ be a saturated group on $\cT$. It is called \emph{locally normal}, if
\begin{enumerate}
\item
the (transitive) subgroup $H(v)$ of $\Sigma_{k(v)}$, which represents the action
of $\St_m(G)$
on the branches coming from a vertex $v$ (supposing that the branching index of $v$
is $k(v)$ and $v$ is at the level $m$), is normal in $\Sigma_{k(v)}$ for any vertex $v\in\cT$;
\item
for any automorphism $\phi:G\to G$ which is defined by an isometry $t$ of $\cT$
fixing $v$, the corresponding element of $\Sigma_{k(v)}$ belongs to $H(v)$.
\end{enumerate}
\end{dfn}

The following statement is evident.

\begin{lem}\label{lem:binarylocnorm}
If $\cT$ is a binary, then $G$ is always locally normal.
\end{lem}

\begin{lem}\label{lem:prosymmgroups}
Let $k\in\N$, $k>1$, $k\ne 4$. Suppose, $H\ss \Sigma_k$ is a normal transitive subgroup
of the symmetric group of permutations on $k$-set. Let $g\in H$ be a non-trivial element.
Then the normal subgroup $N_H(g)$ of $H$ generated by $g$ is transitive.
\end{lem}

\begin{proof}
The group $H$ is $\Sigma_k$ or $A_k$ (the alternating group), while the last one is simple
for these $k$. Hence, $N_H(g)$ is $\Sigma_k$ or $A_k$. It is transitive.
\end{proof}



\begin{lem}\label{lem:stabiliztnaputi}
Suppose, $\phi$ is induced by $t\in \Iso \cT$ and $(v_0,v_1,\dots)$ is an end of
$\cT$. Then for any $n\in\N$ there exists an element $\a_n\in G$ such that
$\a_n t(v_i)=v_i$, $i=1,\dots,n$.
\end{lem}

\begin{proof}
We construct $\a_n$ inductively as a composition $\a_n=\b_{n-1}\dots \b_0$,
$\b_i\in \St_i(G)$. Since $\St_0(G)=G$ acts transitively on the level 1, one can find
$\b_0\in\St_0(G)=G$ such that $\b_0 t(v_1)=v_1$.
Since $\St_1(G)$ acts transitively on the level 2, one can find
$\b_1\in\St_1(G)$ such that $\b_1 \b_0 t(v_2)=v_2$. Moreover, since
$\b_1\in\St_1(G)$, $\b_1 \b_0 t(v_1)=\b_1 v_1=v_1$. And so on.
\end{proof}

\begin{dfn}\label{dfn:weakbranind}
Let $G$ be a weakly branch group. Consider the level $m$. Define
the \emph{weak branch index of the level} $m$ as
$$
\WBI(m)=\max_{v} \WBI(\cT_v),
$$
where $v$ runs over vertexes of the level $m$, while $\WBI(\cT_v)$
is equal to the minimal level in $\cT_v$ of non-trivial action of $G[v]$.
\end{dfn}

\begin{lem}\label{lem:glavdlialocnorm}
Suppose, a saturated weakly branch locally normal group $G$ acts on a
(spherically symmetric)
tree $\cT$ with no vertex of branching index $4$.
Let $t\in\Iso\cT$ induce an automorphism $\phi:G\to G$.
Then for any $m$
there exists an element $\wh g\in\St_m (G)$ such that
$\{\wh g\}_\phi\cap \St_{m+\WBI(m)}(G)=\varnothing$, provided that
$t$ has a fixed vertex at the level $m$.
\end{lem}

\begin{proof}
Let $v_0$ be this fixed vertex.
In the locally normal case
we can make $g't$ act on the
first step successors of $v_0$ without fixed points for some $g'\in \St_m(G)$
(hence, the unique fixed point on $\cT(v_0)$ is $v_0$). If $h g't h^{-1} t^{-1}=g'_h$ is
still in the stabilizer of the next level, then we continue as follows.

By the weakly branch condition one can find a non-trivial element $g''\in G[v]$.
Let $m'\le \WBI(m)$ be its first non-trivial level. Let $v$ be a vertex of the
level $m-1$, such that $g''$ moves its first step successors.
Consider the permutation group $H$ of these successors obtained by the (transitive)
action of $\St_m(G)$. This
group is normal (by the local normality) and contains the representing element of $g''$.
The action of the normalizer of $g''$ in $\St_m(G)$ is a subgroup of $G[v]$ and
its representation on the first step successors of $v$ is transitive by
Lemma \ref{lem:prosymmgroups}. Hence, there is an element $\til g\in G[v]$ such that
$\til g g' t$ has at least one fixed point at the level $m+m'$ while $t$ does not have.
Hence they can not be conjugate by any $h$ and
$(\til g g')_h\not\in \St_{m+\WBI(m)}(G)$ for any $h$. So we are done.
\end{proof}

\begin{teo}\label{teo:mainforlocnormal}
Suppose, a saturated weakly branch locally normal group $G$ acts on a
(spherically symmetric)
tree $\cT$ with no vertex of branching index $4$. Then $R(\phi)=\infty$
for any automorphism of $G$.
\end{teo}

\begin{proof}
Let us take an arbitrary $n\in \N$. We will prove that $R(\phi)\ge n$.
By Lemma \ref{lem:stabiliztnaputi}  find an element $g\in G$, such that
$gt$ has a fixed vertex at each level $1,\dots, w(n)$, where
$$
w(n)=\WBI(1)+\WBI(\WBI(1))+\dots+
\underbrace{\WBI(\WBI(\WBI(\cdots(\WBI}_{n\ {\rm times}}(1))\cdots)))
$$

By Corollary \ref{lem:innerreidem} $R(\phi)=R(\phi')$, where $\phi'$
is induced by $gt$.

We apply inductively $n$ times Lemma \ref{lem:glavdlialocnorm} to prove that
$R(\phi')\ge n$.
\end{proof}

If $\cT$ is a binary tree, then Lemma \ref{lem:binarylocnorm} shows that
$G$ is locally normal. Thus, by Theorem \ref{teo:lavnik}
from Theorem \ref{teo:mainforlocnormal} we obtain the following statement.

\begin{teo}\label{teo:mainbinary}
Let $G$ be a saturated weakly branch group on a binary tree $\cT$. Then
$G$ has the $R_\infty$-property.
\end{teo}

\section{Ternary trees and the Gupta-Sidki group}
The Gupta-Sidki group $G$ \cite{GuSi,SidkiJA1,SidkiJA2} acts on the ternary
tree $\cT$ (see Fig.~\ref{fig:terntree})
\begin{figure}[ht]
\unitlength=1mm
\begin{picture}(100,60)(0,0)
\thicklines
\put(48,50){$\varnothing$}
\put(10,36){\line(3,1){38}}
\put(50,36){\line(0,1){12}}
\put(90,36){\line(-3,1){38}}
\put(9,32){$0$}
\put(49,32){$1$}
\put(89,32){$2$}
\put(0,15){\line(1,2){8}}
\put(10,15){\line(0,1){15}}
\put(20,15){\line(-1,2){8}}
\put(40,15){\line(1,2){8}}
\put(50,15){\line(0,1){15}}
\put(60,15){\line(-1,2){8}}
\put(80,15){\line(1,2){8}}
\put(90,15){\line(0,1){15}}
\put(100,15){\line(-1,2){8}}
\put(-2,10){$00$}
\put(8,10){$01$}
\put(18,10){$02$}
\put(38,10){$10$}
\put(48,10){$11$}
\put(58,10){$12$}
\put(78,10){$20$}
\put(88,10){$21$}
\put(98,10){$22$}
\put(50,5){$\vdots$}
\end{picture}
\caption{}\label{fig:terntree}
\end{figure}
with generators $x$ and $\g$:
$$
x: 0 \to 1 \to 2 \to 0,\quad 00\to 10 \to 20 \to 00,
\dots \quad 0s\to 1s \to 2s \to 0s,
$$
where $s$ is any finite sequence on $0,1,2$ and $\g$ is presented
on Fig.~\ref{fig:gamma}.
\begin{figure}[ht]
\unitlength=1mm
\begin{picture}(45,45)(0,0)
\thicklines
\put(30,35){\line(0,1){10.}}
\put(40,35){\line(-1,1){10}}
\put(29,31){$x$}
\put(39,31){$x^{-1}$}
\put(20,25){\line(0,1){10.}}
\put(30,25){\line(-1,1){10}}
\put(19,21){$x$}
\put(29,21){$x^{-1}$}
\put(0,15){\line(1,1){30}}
\put(10,15){\line(0,1){10.}}
\put(20,15){\line(-1,1){10}}
\put(9,10){$x$}
\put(19,10){$x^{-1}$}
\put(10,5){$\vdots$}
\end{picture}
\caption{}\label{fig:gamma}
\end{figure}

Consider $g=x^{-1}\g^{-1}x\g$. One checks up directly that $g\in \St_1(G)$ and
$g$ acts without fixed points on the level $2$. Now let $(g,g,g)\in \Iso \cT$.
It evidently stabilizes the level $2$ and
acts without fixed points on the level $3$. Moreover, $(g,g,g)\in G$. This is proved
in \cite[Theorem 1]{SidkiJA2}, but also can be seen immediately from the facts that
$(\g,\g,\g)\in G$ and $(x,x,x)\in \Iso\cT$ is in the normalizer of $G$. Hence,
$(g,g,g)\in \St_2(G)$. By \cite[Theorem 1]{SidkiJA2} (note that there is a misprint
in the formulation of that theorem: $i+1$ should be replaced by $i-1$ three times)
$g_i:=(g,g,\dots,g) \in \St_i(G)$ for the corresponding $i$.
In particular, the subgroup of $\Sigma_3$ mentioned in the first item of
Definition \ref{dfn:normality} contains $A_3\cong \Z_3$. On the other hand,
$G$ is constructed by actions $1$, $x$, and $x^{-1}$, hence, it is not larger
than $A_3$. Thus the first item of Definition \ref{dfn:normality} holds.

Let us remind the description of automorphisms of $G$ \cite{SidkiJA2} (Theorem 3
and pp. 39--41): $\Aut(G)=(G \rtimes X)\rtimes V$, where $X$ is an elementary
abelian $3$-group of infinite rank with basis $x^{(i)}$ ($i\in \N$):
$$
x^{(1)}=(x,x,x),\quad x^{(i+1)}=(x^{i},x^{i},x^{i}) \quad {\rm for}\ i\ge 1,
$$
and $V\cong \Z_4$ with nontrivial elements $\t_1$, $\t_2$, $\t_3$, which
map generators of $G$ in the following way
$$
\begin{array}{ll}
\t_1(\g)=\g, & \t_1(x)=x^{-1},\\
\t_2(\g)=\g^{-1},& \t_2(x)=x,\\
\t_3(\g)=\g^{-1}, & \t_3(x)=x^{-1},
\end{array}
$$
and act on the base of $X$ by the formulas ($\t_3=\t_2\circ\t_1$):
$$
\t_1(x^{i})=(x^{i})^{-1},\quad
\t_2(x^{2i-1})=(x^{2i-1})^{-1},\quad \t_2(x^{2i})=(x^{2i}),
$$
for $i\ge 1$. Thus the action of $t$, mentioned in the second item of
Definition \ref{dfn:normality}, can be only defined by $1$, $x$ or $x^{-1}$,
i.e., belongs to $A_3$. So, $G$ is locally normal.

Also, this group is saturated weakly branch \cite[Prop. 8.6]{LavrNekr},
\cite{SidkiJA2}. Hence applying Theorem \ref{teo:mainforlocnormal} we obtain
\begin{teo}\label{teo:guptasidki}
For any automorphism $\phi$ of the Gupta-Sidki group one has $R(\phi)=\infty$.
\end{teo}

\begin{rk}\label{rk:GupSidstrsat}
The argument above concerning the first property of Definition
{\rm\ref{dfn:normality}} shows that the Gupta-Sidki group is strongly
saturated (Definition \ref{dfn:strongsatur}).
\end{rk}

In fact, the unique proper transitive subgroup of $\Sigma_3$ is $A_3$. Hence, one has
the following statement.

\begin{teo}
Any saturated group on a ternary tree enjoys the first property of Definition
{\rm\ref{dfn:normality}}.
\end{teo}

\def\cprime{$'$} \def\dbar{\leavevmode\hbox to 0pt{\hskip.2ex \accent"16\hss}d}
  \def\polhk#1{\setbox0=\hbox{#1}{\ooalign{\hidewidth
  \lower1.5ex\hbox{`}\hidewidth\crcr\unhbox0}}}
\providecommand{\bysame}{\leavevmode\hbox to3em{\hrulefill}\thinspace}
\providecommand{\MR}{\relax\ifhmode\unskip\space\fi MR }
\providecommand{\MRhref}[2]{%
  \href{http://www.ams.org/mathscinet-getitem?mr=#1}{#2}
}
\providecommand{\href}[2]{#2}

\end{document}